\documentclass[10pt]{amsart}

\usepackage{amsmath, amsfonts, amssymb, xypic}
\usepackage{latexsym, verbatim}
\usepackage{graphicx}
\usepackage{color}
\usepackage{url}
\usepackage{soul}
\usepackage[table]{xcolor}
\usepackage{hyperref}
\usepackage{amsmath}

\DeclareMathOperator{\End}{End}
\DeclareMathOperator{\csch}{csch}

\theoremstyle{definition}

\newcommand{\N}{\mathbb{N}}
\newcommand{\Z}{\mathbb{Z}}
\newcommand{\R}{\mathbb{R}}

\def\Aut{\operatorname{Aut}}

\title{Almost complex manifolds with small Nijenhuis tensor}

\author[L. Fernandez]{Luis Fernandez}
  \address[L. Fernandez]{Department of Mathematics and Computer Science, Bronx Community College, City University of New York, 2155 University Ave., Bronx, NY 10453}
  \email{luis.fernandez01@bcc.cuny.edu}

\author[T. Shin]{Tobias Shin}
  \address[T. Shin]{Department of Mathematics, Stony Brook, State University of New York, 100 Nicolls Rd., Stony Brook, NY 11794}
  \email{tobias.shin@stonybrook.edu}
   
  \author[S. Wilson]{Scott O. Wilson}
  \address[S. Wilson]{Department of Mathematics, Queens College, City University of New York, 65-30 Kissena Blvd., Flushing, NY 11367}
  \email{scott.wilson@qc.cuny.edu}
  
\keywords{almost complex manifold, Nijenhuis tensor}
\subjclass[2010]{32Q60, 53C15}

\begin{document}

\begin{abstract} We give several explicit examples of compact manifolds with a $1$-parameter family of almost complex structures having arbitrarily small Nijenhuis tensor in the $C^0$-norm. The $4$-dimensional examples possess no complex structure, whereas the $6$-dimensional example does not possess a left invariant complex structure, and whether it possesses a complex structure appears to be unknown.
\end{abstract}

\maketitle

%\allowdisplaybreaks
%\setcounter{tocdepth}{2}
%\setcounter{tocdepth}{1}
%\tableofcontents

\section{Introduction}

The purpose of the following short note is to give several explicit examples of compact manifolds with a $1$-parameter family of almost complex structures having arbitrarily small Nijenhuis tensor in the $C^0$, or supremum, norm. The $4$-dimensional examples possess no complex structure, whereas the $6$-dimensional example does not possess a left invariant complex structure, and whether it possesses a complex structure appears to be unknown.

The idea that such examples might exist was inspired by efforts of the second author to place the work of Demailly and Gaussier \cite{DG} in the context of Gromov's h-principle, whereby an integrable complex structure was interpreted as a holonomic solution of a locally closed differential relation.  The lack of homotopy obstructions to formal solutions of this differential relation led the second author to attempt various h-principle techniques, to try to deform a ``formal integrable complex structure'' into a genuine one. This naturally led to the question of whether every almost complex manifold has almost complex structures that are arbitrarily close to an integrable one.

We remark  that the examples below appear not to be isolated, as some exist in larger parameter families than we've presented here. The technique used to construct these examples involves considerable trial-and-error, aided by computer algebra software, and in some cases guided by the gradient descent method and the assumption of rational functions. It remains an interesting practical problem to find a general technique. 

Perhaps even more importantly, one wishes 
 to have a more clear conceptual reason as to why this is possible, and it remains an open question whether almost complex structures with arbitrarily $C^0$-small Nijenhuis tensor always exist on compact almost complex manifolds. Some recent results in this direction, particular to dimension $6$, appear in the work of Fei et. al. \cite{FPPZ}, establishing sufficient conditions for some symplectic manifolds.

\medskip
\noindent
\textbf{Acknowledgments.} The authors would like to thank Paul Feehan for encouraging us to share these examples in the present paper, and to thank the referee for making several corrections and improvements.

\section{Nilmanifold examples}

We begin with two nilmanifold examples, in real dimensions $4$ and $6$, respectively. In each case, the nilpotent Lie algebra $\mathfrak g$ has rational structure constants, so there exists a lattice $\Gamma$ of the Lie group $G$ that has compact 
quotient $\Gamma  \backslash G$ (see \cite{Malcev}, Theorem 7). Then any linear complex structure on $\mathfrak g$ descends to a left invariant almost complex structure on any such quotient $\Gamma  \backslash G$. Any two norms on the finite dimensional vector space $\mathfrak g$ are equivalent, and convergence of a tensor will be understood in the induced $C^0$-topology, which is norm-independent.

In general, the integral cohomology of  $\Gamma  \backslash G$ depends on the lattice $\Gamma$, but the real cohomology does not, and can be computed from the cohomology of the Lie algebra, by Nomizu's theorem \cite{Nomizu}. The diffeomorphism type of the quotient $\Gamma \backslash G$ is completely determined by the fundamental group, i.e. the lattice $\Gamma$. In fact, this is true for solvmanifolds as well, see \cite{Mos}.

\subsection{Filiform Lie algebra in Dimension $4$}\label{4dfiliform}

Consider the real $4$-dimensional nilpotent Lie algebra $\mathfrak g$, with basis $\{X_1 , X_2, X_3, X_4\}$, and only non-zero brackets determined by
\[
[X_1 ,X_i  ] = X_{i+1} \quad \textrm{for} \quad i=2,3,
\]
or similarly, if $\{x_1,x_2,x_3,x_4\}$ is the dual basis,
\[
dx_1=0 \quad  dx_2 = 0 \quad dx_3 = -x_1\wedge x_2 \quad dx_4 = -x_1\wedge x_3. 
\]

A \emph{filiform $4$-manifold} $M^4$ is a compact quotient $ \Gamma  \backslash G$, where $\Gamma$ is a lattice of  the simply connected Lie group $G$  of $\mathfrak g$.  Using the definition of $\mathfrak g$ above, it is easy to see that $x_1$ and $x_2$ generate $H^1(M^4)$, and $x_1\wedge x_4$ and $x_2\wedge x_3$ generate $H^2(M^4)$, so that the Betti numbers of $M^4$ are $b_1= b_2 = b_3 = 2$.

The manifold $M^4$ does not admit any integrable complex structure. Indeed, since $b_1$ is even, by Kodaira's 
classification of surfaces, $M^4$ would then be K\"{a}hler, and hence $M^4$ would be formal. But this manifold is not formal, as it has a nontrivial Massey product, and moreover, every formal nilmanifold is diffeomorphic to a torus \cite{Ha89}, and thus has $b^1=4$.  Alternatively, one can argue that since $b_1=2$ and $M^4$ is parallelizable, then by a result of Fern\'andez and Gray \cite{FG} (which also relies on the classification of surfaces)
$M^4$ does not have a complex structure.  We note that this manifold does admit a symplectic form. One example is $\omega = x_1 \wedge x_4 + x_2 \wedge x_3$. 

We give an example of a $1$-parameter family $J_t$ of left-invariant almost-complex structures on $M^4$ such that
the Nijenhuis tensor $N_t:=N(J_t)$ satisfies $N_t\to 0$ as $t \to \infty$. In the ordered basis $\{X_1 , X_2, X_3, X_4\}$, define
\[
J_t =
\begin{bmatrix}
1& -2 \csch t& 0& 0 \\ 
\sinh t& -1& 0& 0 \\
0& 0& -1 - \sqrt 2 & -2 (2 + \sqrt 2 ) \csch t \\
0& 0& \sinh t & 1 + \sqrt 2 \\
\end{bmatrix},
\]
where $\csch t = 1/ \sinh t$.

Since the Nijenhuis tensor 
\[
N_t(X,Y) := [J_tX,J_tY]  - [X,Y] - J_t[X,J_tY] - J_t [J_tX,Y]
% [X,Y]+J_t ([J_tX,Y]+J_t[X,J_tY])-[J_tX,J_tY]
\]
is skew-symmetric, it is completely determined by

\begin{align*}
N_t( X _1, X _2) & =0 \\
N_t( X _1, X _3) & =  -4(1+ \sqrt 2 ) \csch  t \; X _3 \\
N_t( X _1, X _4) & = -4\csch  t \, \left( 2 (3 + 2\sqrt 2) \csch  t \; X _3 - (1+ \sqrt 2) \; X _4 \right) \\
N_t( X _2, X _3) &= -4 \csch  t\, \left((2 + \sqrt 2) \csch  t \; X _3 - (1 + \sqrt 2) \; X _4 \right) \\
N_t( X _2, X _4) &= 4 (2 + \sqrt 2) \csch^2  t  \; X _4 \\
N_t( X _3, X _4) & = 0.
\end{align*}
It is clear that every component approaches zero uniformly as $t \to \infty$, and therefore $N_t \to 0$.

Thus we see any such filiform $4$-manifold $\Gamma \backslash G$ has an almost K\"ahler structure, and has an arbitrarily small Nijenhuis tensor for another $J$, but no complex structure. 

\subsection{Filiform Lie algebra in Dimension $6$}\label{6dfiliform}

Consider the filiform Lie algebra $\mathfrak g$ of dimension $6$, with basis $\{X_1 , \ldots ,X_6\}$, and only non-zero brackets determined by
\[
[X_1 ,X_i  ] = X_{i+1} \quad \textrm{for} \quad i=2,3,4,5.
\]
A  \emph{filiform $6$-manifold} $M^6$ is a compact quotient $\Gamma  \backslash G$ 
where $\Gamma$ is a lattice in the simply connected Lie group $G$ associated to $\mathfrak g$. 

According to \cite{GM}, the Lie algebra  $\mathfrak g$ does not admit an integrable linear complex structure, and thus there is no left invariant complex structure on any of the compact quotients $M^6$. To our knowledge, there is no known complex structure on these manifolds, nor any proof that none of these compact 
quotients admit a complex structure (as this is a widely open problem for all almost complex $6$-manifolds). Incidentally, such manifolds do admit an almost K\"ahler structure, e.g. 
\[
JX_1= X_6, \quad JX_2=X_5, \quad JX_3=-X_4, 
\]
with 
\[
\omega = x_1 \wedge x_6 + \wedge x_2 x_5 - x_3 \wedge x_4,
\] 
where $\{x_1, \dots , x_6\}$ is the dual basis.

We give an example of a $1$-parameter family $J_t$, of left-invariant almost-complex structures on any $M^6$, such that
the Nijenhuis tensor $N_t:=N(J_t)$ satisfies $N_t\to 0$ as $t \to \infty$. In the ordered basis $\{X_1 , \ldots, X_6\}$, define
$$\arraycolsep=2.5pt\renewcommand*{\arraystretch}{1.5}
J_t = \left[
\begin{array}{cccccc}
 \frac{t^2-1}{\sqrt{3} \left(t^2+1\right)} & -\frac{1}{t^3} & 0 & 0 & 0 & 0 \\[5pt]
 \frac{4 t^3\left(t^4+t^2+1\right)}{3 \left(t^2+1\right)^2} & -\frac{t^2-1}{\sqrt{3} \left(t^2+1\right)} & 0 &
   0 & 0 & 0 \\
 1 & 0 & \frac{1}{\sqrt{3}} & \frac{1}{t} & 0 & 0 \\
 -\frac{2 t^3}{\sqrt{3} \left(t^2+1\right)} & \frac{1}{t^2} & -\frac{4 t}{3} & -\frac{1}{\sqrt{3}} & 0 & 0 \\[2pt]
 0 & -\frac{5 t^4+4 t^2+1}{\sqrt{3} t^3 \left(2 t^2+1\right)} & 
\frac{8 t^2 (t^2+1)}{3 \sqrt{3} (2 t^2+1)}&
 -\frac{4 t^5}{6 t^4+9 t^2+3} & \frac{1}{\sqrt{3}} & \frac{t^2+1}{t^3} \\[5pt]
 \frac{4 t^3\left(t^4+t^2+1\right)}{3 \sqrt{3} \left(2 t^4+3 t^2+1\right)} & \frac{2}{3} &
-\frac{16 t^5 \left(2 t^4+2 t^2+1\right)}{9 \left(t^2+1\right)^2 \left(2 t^2+1\right)}
& -\frac{8 t^4}{\sqrt{3} \left(6
   t^2+3\right)} & -\frac{4 t^3}{3 \left(t^2+1\right)} & -\frac{1}{\sqrt{3}}
\end{array}
\right].
$$
Then it is a long but straightforward algebraic computation to find the entries of the tensor $N_t$ as given below:
\begin{eqnarray*}
N_t (X_1,X_2)&=&\frac{1}{t^3}X_4+\frac{\left(t^6+7 t^4+5 t^2+1\right)}{\sqrt{3} t^6 \left(2
   t^2+1\right)}  X_5-\frac{2 \left(t^2+1\right)}{3 t \left(2
   t^2+1\right)} X_6\\
N_t (X_1,X_3)&=&-\frac{2 t}{\sqrt{3} \left(t^2+1\right)} X_3-\frac{4}{3
   \left(t^2+1\right)} X_4 - \frac{8 \left(t^4+t^2+1\right)}{3 \sqrt{3} t
   \left(t^2+1\right)^2} X_5\\
N_t (X_1,X_4)&=&-\frac{1}{t^2}X_3+\frac{2 t}{\sqrt{3}
   \left(t^2+1\right)} X_4\\
N_t (X_1,X_5)&=&-\frac{2
   }{\sqrt{3} t} X_5-\frac{4}{3 \left(t^2+1\right)} X_6\\
N_t (X_1,X_6)&=&-\frac{\left(t^2+1\right)^2
   }{t^6}X_5+\frac{2}{\sqrt{3} t} X_6\\
N_t (X_2,X_3)&=&\frac{1}{t^4}X_3-\frac{2}{\sqrt{3} t^3} X_4+\frac{4 \left(t^4+3
   t^2+1\right)}{3 t^2 \left(t^2+1\right) \left(2 t^2+1\right)} X_5-\frac{8}{3 \sqrt{3}
   t} X_6\\
N_t (X_2,X_4)&=&-\frac{1}{t^4}X_4+\frac{2}{\sqrt{3} t^3} X_5-\frac{4}{3 \left(2
   t^2+1\right)} X_6\\
N_t (X_2,X_5)&=&\frac{\left(t^2+1\right)}{t^6} X_5-\frac{2}{\sqrt{3}
   t^3} X_6\\
N_t (X_2,X_6)&=&-\frac{\left(t^2+1\right)}{t^6}X_6\\
N_t (X_i,X_j)&=&0 \ \ \mbox{if $i,j\ge 3$},
\end{eqnarray*}
and so $N_t \to 0$ as $t \to \infty$.

\section{Solvmanifold examples}

We give two families of examples of solvmanifolds of real dimension $4$ which have no complex structures yet have a family of almost complex structures whose Nijenhuis tensors tend to zero. 
Since there is no general criterion to ensure the existence of a co-compact lattice in a simply connected solvable Lie group, we will specify one which is co-compact in each case. As before, since the complex structures on $\mathfrak g$ that we provide are linear, they descend to any compact quotient.

\subsection{A first class of $4d$-solvmanifold examples}

The first solvmanifold is taken from Fern\'andez and Gray \cite{FG}.  Consider the real $4$-dimensional  solvable Lie algebra $\mathfrak g$, with basis $X_1 , X_2, X_3, X_4$, and only non-zero brackets
\begin{align*}
[X_1 ,X_3  ] = -k X_{1} \\
[X_2 ,X_3  ] = k X_{2} ,
\end{align*}
for any $k \neq 0$.    This Lie algebra is the direct sum of the trivial one dimensional Lie algebra (generated by $X_4$) with the Lie algebra $\mathfrak g(k)$ of the simply connected solvable (non-nilpotent) Lie group $G(k)$ given by matrices of the form
\[
\begin{bmatrix}
e^{kz} & 0 & 0 & x \\
0 & e^{-kz} & 0 & y \\
0 & 0 & 1 & z \\
0 & 0 & 0 & 1
\end{bmatrix}.
\]
For this case, we take $x,y,z \in \R$ and $k$ a real number such that $e^{k} + e^{-k} $ is an integer different than $2$ (so that
$k \neq 0$).

Let $M^4(k) = (\Gamma \backslash G(k)) \times S^1$, where $\Gamma$ 
is any lattice of $G(k)$ determined, as in (\cite{AGH}, Theorem 4 (4)), by the subgroup of $G(k)$ generated by matrices of the form
\[
\begin{bmatrix}
1 & 0 & 0 & u_1 \\
0 & 1 & 0 & u_2 \\
0 & 0 & 1 & 0 \\
0 & 0 & 0 & 1
\end{bmatrix},
\quad
\begin{bmatrix}
1 & 0 & 0 & v_1 \\
0 & 1 & 0 & v_2\\
0 & 0 & 1 & 0 \\
0 & 0 & 0 & 1
\end{bmatrix},
\quad
\begin{bmatrix}
e^{kn} & 0 & 0 & 0 \\
0 & e^{-kn} & 0 & 0 \\
0 & 0 & 1 & n \\
0 & 0 & 0 & 1
\end{bmatrix},
\]
where $n \in \N$, and $(u_1,u_2)$ and $(v_1,v_2)$ are linearly independent. 

This Lie algebra is completely solvable, i.e. the adjoint action $ad_X$ has real eigenvalues for all $X$.
By a theorem of Hattori \cite{Hatt}, one can compute the de Rham cohomology of any associated completely solvmanifold $\Gamma \backslash G$ from the cohomology of the Lie algebra. 

Using the brackets above, it is easy to see that $x_3$ and $x_4$ generate $H^1(M^4)$, and $x_1\wedge x_2$ and $x_3\wedge x_4$ generate $H^2(M^4)$, so that the Betti numbers of $M^4$ are $b_1= b_2 = b_3 = 2$.
Since $b_1=2$ and $M^4(k)$ is parallelizable, then by a result of Fern\'andez and Gray \cite{FG}, $M^4(k)$ does not have a complex structure. As also pointed out in \cite{FG}, these manifolds nevertheless are formal, symplectic,  and moreover satisfy all the known cohomological  properties of a K\"ahler manifold.

Since the manifolds $M^4(k)$ all have the same minimal model as $S^1 \times S^1  \times S^2$, as shown in \cite{FG}, this shows that there is no algebraic condition on the minimal model of four manifolds implying the existence of a complex structure. A. Milivojevic obtained a similar result in dimension $2n=6$ and greater using a geometric argument, namely, for every almost complex manifold, he constructs a non-almost-complex manifold by taking its connected sum with a non-spin$^c$ simply connected rational homology sphere  \cite{M}.

For any $k \neq 0$, consider the family of linear almost complex structures on  $\mathfrak g$ defined 
in the ordered basis $\{X_1 , \ldots, X_4\}$ by 
\[
\renewcommand\arraystretch{2}
J_t =
\begin{bmatrix}
\frac{-2 }{kt^2} & \frac{-1}{\sqrt3} &  -\frac{6+\sqrt3 k t^2+2k^2 t^4}{3k^2 t^3} &  \frac{6-\sqrt3 k t^2+2k^2 t^4}
{3k^2 t^5}  \\
\frac{-1}{\sqrt 3} & 0 & \frac{-1}{\sqrt 3 k t} -\frac{2t}{3} &  \frac{\sqrt 3 -2kt^2}{3kt^3} \\
\frac{1}{t} &  \frac{1}{t} & \frac{1}{\sqrt 3} + \frac{1}{kt^2} &  \frac{-1}{kt^4} \\
-t &  t & \frac{-1}{k} &  \frac{-1}{\sqrt 3} + \frac{1}{kt^2}  \\
\end{bmatrix}.
\]

Then the Nijenhuis tensor is determined by 
\begin{align*}
N_t( X _1, X _2) & =\frac{-2k}{\sqrt 3t}\left( X_1 +  X_2 \right) + \frac{2k}{t^2}X_3  \\
N_t( X _1, X _3) & = \frac 1 t \left(\frac{2k}{\sqrt 3}X_3 + X_4\right) +\frac{1}{t^2}\left(\frac{-1}{\sqrt 3}X_1
+\frac{1}{\sqrt 3}X_2\right) 
-\frac{1}{t^3}X_3 + \frac{2}{kt^4}X_1  \\
N_t( X _1, X _4) & = \frac{2}{kt^6}X_1 -\frac{1}{t^5}X_3+\frac{1}{\sqrt 3 t^4}\left(X_2-X_1\right)
+ \frac{1}{t^3}X_4 + \frac{2k}{3t^2}\left(X_1+X_2\right)
\\
N_t( X _2, X _3) & = \frac{-2}{kt^4}X_1 -\frac{1}{t^3}X_3+\frac{1}{\sqrt3 t^2}\left(3X_1+X_2\right)
+\frac{1}{t}\left(\frac{-2k}{\sqrt3}X_3-X_4\right) \\
N_t( X _2, X _4) & = \frac{2}{kt^6}X_1+\frac{1}{t^5}X_3-\frac{1}{\sqrt3 t^4}\left(3X_1+X_2\right)+\frac{1}{t^3}X_4+\frac{2k}{3t^2}\left(X_1+X_2\right)\\
N_t( X _3, X _4) & = \frac{4}{k^2t^7}X_1-\frac{2}{kt^6}X_3+\frac{2}{\sqrt3 k t^5}\left(X_2-X_1\right)
+\frac{2}{t^4}\left(\frac{1}{\sqrt 3}X_3+\frac 1 k X_4\right) \\
& + \frac{2}{3t^3}\left(X_1+X_2\right)
-\frac{4k}{3t^2}X_3 +\frac{4 k}{3\sqrt 3 t}\left(X_1+X_2\right).
\end{align*}
 Then, for each $k \neq 0$, $N_t \to 0$ as $t \to \infty$, as claimed.  

\subsection{A second class of $4d$-solvmanifold examples}

The next class of almost complex manifolds that admit no complex structure are taken from Hasegawa (\cite{Ha05}, Section 4, Example 2). Consider the Lie group $\R^3 \rtimes_\phi \R$, where
$\phi: \R \to \Aut(\R^3)$, and $\R$ and $\R^3$ have the standard additive structures. In order to ensure there is a lattice, we choose $\phi(1) = A \in SL(3,\Z)$, and restrict to the case that
$A$ has three real, positive, distinct eigenvalues. There are many such examples, e.g. 
\[
A= 
\begin{bmatrix}
0 & 0 & 1 \\
1 & 0 & -k \\
0 & 1 & 8
\end{bmatrix}
\]
for any integer $k$ with $6 \leq k \leq 15$, (c.f. Bock \cite{Bock}). Then the action of $A$ on $\Z^3$ can be extended to  $\phi: \R \to \Aut(\R^3)$ by defining $\phi_t= \exp( t \log A)$, and 
$\Gamma = \Z^3 \rtimes_\phi \Z$ is a lattice in $G= \R^3 \rtimes_\phi \R$. Since the eigenvalues of $A$ are distinct, we can 
write $A = VDV^{-1}$, with $D$ diagonal having entries $e^{\lambda_1},e^{\lambda_2}, e^{-(\lambda_1+\lambda_2)}$,
and $V$ is a matrix whose columns are the respective eigenvectors.  Let $X_1$ be the standard basis vector for $\R$, and let $\{X_2, X_3, X_4\}$
be the columns of $V^{-1}$.

The derivative of $\phi: \R \to \Aut(\R^3)$ at zero is the map $\phi: \R \to \End(\R^3)$ sending $X_1$ to left multiplication by $\log A$.  So, in the standard basis $\{E_i\}$ of $\R \times \R^3$, with $E_1 = X_1$, the Lie algebra of $G$ has non-zero brackets determined by 
\[
[E_1,E_i] = \left( \log A \right)  E_i,
\]
for $i=2,3,4$. Equivalently, the 
transport of the Lie bracket by $V$, defined by $\{ \, , \,\} : = V^{-1} \circ [ V (-), V(-) ]$, satisfies 
\begin{align*}
\{X_1 ,X_2  \} &= \lambda_1 X_{2} \\
\{X_1 ,X_3  \}&= \lambda_2 X_{3} \\
\{X_1 ,X_4  \} &= -(\lambda_1 + \lambda_2) X_{4}.
\end{align*}
 Note that the first Betti number of $\Gamma \backslash G$ is equal to one unless any of $\lambda_1, \lambda_2,$ or
 $(\lambda_1 + \lambda_2)$ are zero. 
 
 According to Hasegawa's classification of compact complex $4$-dimensional solvmanifolds, any such
 solvmanifold $\Gamma \backslash G$ does not have a complex structure (\cite{Ha05}, Section 4, Example 2).
 
 To give a $1$-parameter family of almost complex structures $J_t$ on $\left(\mathfrak{g}, [ \,, \,]\right)$ with 
 $N_{J_t} \to 0$, it suffices to give almost complex structures $K_t$ on $\left(\mathfrak{g}, \{ \,, \, \} \right)$ 
 in the basis $\{X_1, X_2, X_3, X_4\}$, with $N_{K_t} \to 0$, for then we may define $J_t:=V^{-1} K_t V$. To this end, in the ordered basis $\{X_1 , \ldots, X_4\}$, let
\[
\renewcommand\arraystretch{2}
K_t =
\begin{bmatrix}
1& 1/t & \frac{2(\lambda_1+2\lambda_2)}{(\lambda_1-\lambda_2)t}& 0 \\ 
-\frac{2(2\lambda_1+\lambda_2)t}{\lambda_1-\lambda_2} & -\frac{2\lambda_1+\lambda_2}{\lambda_1-\lambda_2}& -\frac{2(2\lambda_1+\lambda_2)(\lambda_1+2\lambda_2)}{(\lambda_1-\lambda_2)^2} & \frac{\lambda_1+2\lambda_2}{(\lambda_1-\lambda_2)t} \\ 
t& 1/2& \frac{\lambda_1+2\lambda_2}{\lambda_1-\lambda_2}  & \frac{-1}{2t} \\
0 & t& \frac{2(2\lambda_1+\lambda_2)t}{\lambda_1-\lambda_2}& 0 \\ 
\end{bmatrix}
\]

Then the Nijenhuis tensor of $K_t$ is determined by 
\begin{align*}
N_t( X _1, X _2) & = \frac{ \lambda_1+2 \lambda_2 }{t} X_1 \\
N_t( X _1, X _3) & = \frac{2(2\lambda_1+\lambda_2)(\lambda_1+2\lambda_2)}{(\lambda_1- \lambda_2)t}X_1 \\
N_t( X _1, X _4) & = -\left( \frac{ \lambda_1+2 \lambda_2 }{t^2} \right) X_1 
+
 \left(\frac{2(2\lambda_1+\lambda_2)(\lambda_1+2\lambda_2)}{(\lambda_1- \lambda_2)t } \right) X_2 
 - \left( \frac{ \lambda_1+2 \lambda_2 }{t}  \right) X_3 \\
N_t( X _2, X _3) & = \left( \frac{  2\left( \lambda_1+2 \lambda_2 \right) }{t^2}  \right) X_1 
-
 \left(  \frac{2(2\lambda_1+\lambda_2)(\lambda_1+2\lambda_2)}{(\lambda_1- \lambda_2)t} \right) X_2 
 + \left(  \frac{ \lambda_1+2 \lambda_2 }{t} \right) X_3
   \\
N_t( X _2, X _4) & =    \left( \frac{(2\lambda_1+\lambda_2)(\lambda_1+2\lambda_2)}{(\lambda_1- \lambda_2)t^2} \right) X_2  - \left(  \frac{ \lambda_1+2 \lambda_2 }{2t^2} \right)X_3 \\
N_t( X _3, X _4) & = \left( 2 \frac{(2\lambda_1+\lambda_2)(\lambda_1+2\lambda_2)^2}{(\lambda_1- \lambda_2)^2 t^2} \right)X_2 - \left(  \frac{ \left( \lambda_1+2 \lambda_2 \right)^2}{  \left( \lambda_1- \lambda_2 \right)  t^2}  \right)X_3.
\end{align*}
So, $N_t \to 0$ as $t \to \infty$.

\bibliographystyle{alpha}

\bibliography{biblio}

\begin{thebibliography}{FPPZ20}

\bibitem[AGH61]{AGH}
L.~Auslander, L.~Green, and F.~Hahn.
\newblock Flows on some three-dimensional homogeneous spaces.
\newblock {\em Bull. Amer. Math. Soc.}, 67:494--497, 1961.

\bibitem[Boc16]{Bock}
C.~Bock.
\newblock On low-dimensional solvmanifolds.
\newblock {\em Asian J. Math.}, 20(2):199--262, 2016.

\bibitem[DG17]{DG}
J.P. Demailly and H.~Gaussier.
\newblock Algebraic embeddings of smooth almost complex structures.
\newblock {\em J. Eur. Math. Soc. (JEMS)}, 19(11):3391--3419, 2017.

\bibitem[FG90]{FG}
M.~Fern\'{a}ndez and A.~Gray.
\newblock Compact symplectic solvmanifolds not admitting complex structures.
\newblock {\em Geom. Dedicata}, 34(3):295--299, 1990.

\bibitem[FPPZ20]{FPPZ}
T.~Fei, D.~Phong, S.~Picard, and X.~Zhang.
\newblock Geometric flows for the type {IIA} string.
\newblock {\em Preprint arxiv:2011.03662v1}, 2020.

\bibitem[GR02]{GM}
M.~Goze and E.~Remm.
\newblock Non existence of complex structures on filiform {L}ie algebras.
\newblock {\em Comm. Algebra}, 30(8):3777--3788, 2002.

\bibitem[Has89]{Ha89}
K.~Hasegawa.
\newblock Minimal models of nilmanifolds.
\newblock {\em Proc. Amer. Math. Soc.}, 106(1):65--71, 1989.

\bibitem[Has05]{Ha05}
K.~Hasegawa.
\newblock Complex and {K}\"{a}hler structures on compact solvmanifolds.
\newblock {\em J. Symplectic Geom.}, 3(4):749--767, 2005.
\newblock Conference on Symplectic Topology.

\bibitem[Hat60]{Hatt}
A.~Hattori.
\newblock Spectral sequence in the de {R}ham cohomology of fibre bundles.
\newblock {\em J. Fac. Sci. Univ. Tokyo Sect. I}, 8:289--331 (1960), 1960.

\bibitem[Mal49]{Malcev}
A.~I. Mal'cev.
\newblock On a class of homogeneous spaces.
\newblock {\em Izvestiya Akad. Nauk. SSSR. Ser. Mat. Amer. Math. Soc. Transl.,
  39 (1951)}, 13:9--32, 1949.

\bibitem[Mil21]{M}
Aleksandar Milivojevi\'{c}.
\newblock On the realization of symplectic algebras and rational homotopy types
  by closed symplectic manifolds.
\newblock {\em Proc. Amer. Math. Soc.}, 149(5):2257--2263, 2021.

\bibitem[Mos54]{Mos}
G.~D. Mostow.
\newblock Factor spaces of solvable groups.
\newblock {\em Ann. of Math. (2)}, 60:1--27, 1954.

\bibitem[Nom54]{Nomizu}
K.~Nomizu.
\newblock On the cohomology of compact homogeneous spaces of nilpotent {L}ie
  groups.
\newblock {\em Ann. of Math. (2)}, 59:531--538, 1954.

\end{thebibliography}

\end{document}